\documentstyle[fullpage,12pt,psfig]{article}
\def\SBIMSMark#1#2#3{
 \font\SBF=cmss10 at 10 true pt
 \font\SBI=cmssi10 at 10 true pt
 \setbox0=\hbox{\SBF Stony Brook IMS Preprint \##1}
 \setbox2=\hbox to \wd0{\hfil \SBI #2}
 \setbox4=\hbox to \wd0{\hfil \SBI #3}
 \setbox6=\hbox to \wd0{\hss
             \vbox{\hsize=\wd0 \parskip=0pt \baselineskip=10 true pt
                   \copy0 \break%
                   \copy2 \break%
                   \copy4 \break}}
 \dimen0=\ht6   \advance\dimen0 by \vsize \advance\dimen0 by 8 true pt
                \advance\dimen0 by -\pagetotal
 \dimen2=\hsize \advance\dimen2 by .25 true in
  \openin2=publishd.tex
  \ifeof2\setbox0=\hbox to 0pt{}
  \else 
     \setbox0=\hbox to 3.1 true in{
                \vbox to \ht6{\hsize=3 true in \parskip=0pt  \noindent  
                \input publishd.tex 
                \vfill}}
  \fi
  \closein2
  \ht0=0pt \dp0=0pt
 \ht6=0pt \dp6=0pt
 \setbox8=\vbox to \dimen0{\vfill \hbox to \dimen2{\copy0 \hss \copy6}}
 \ht8=0pt \dp8=0pt \wd8=0pt
 \copy8
 \message{*** Stony Brook IMS Preprint #1, #2 ***}
}

\title{Teichm\"{u}ller distance for some polynomial-like maps}
\author{Eduardo A. Prado \thanks{Supported in part by CNPq-Brazil}\\  Instituto de Matem\'atica e Estat\'istica \\
Universidade de S\~ao Paulo \\Caixa Postal 66281 CEP 05389-970 \\
S\~ao Paulo Brazil \\ e.mail: prado@ime.usp.br} 

\date{ }
\newcommand{\apri}{\it a \: priori}

\newcommand{\och}{off-critically hyperbolic}
\newcommand{\ochs}{off-critically hyperbolic }
\newcommand{\diam}{ {\rm diam}}

\newcommand{\julf}{ {\em J(f)}}
\newcommand{\julg}{ {\em J(g)}}
\newcommand{\finjulf}{ {\em K(f)} }
\newcommand{\teichdist}{ {\rm d_{T} } }
\newcommand{\orb}{ {\cal O} }
\newcommand{\pollikef}{f: \bigcup U_{i} \rightarrow U}
\newcommand{\pollikeg}{g: \bigcup V_{i} \rightarrow V}

\newcommand{\preas}{ {\rm P} }
\newcommand{\hausdim}{ {\rm HD} }
\newcommand{\entr}{ {\rm h} }

\newtheorem{dummylemma}{Dummy}[section]
\newtheorem{theorem}[dummylemma]{Theorem}
\newtheorem{thm}[enumi]{Theorem}
\newtheorem{lemma}[dummylemma]{Lemma}
\newtheorem{corollary}[dummylemma]{Corollary}
\newtheorem{ExamplE}[dummylemma]{\em Example}
\newtheorem{DefinitioN}[dummylemma]{\em Definition}

\newenvironment{romanEnum}
                {\begin{list}
                    {(\roman{enumi})}           
                    {\usecounter{enumi}         
                     \setlength{\rightmargin}{\leftmargin}
                     \setlength{\labelwidth}{7.5mm}
                    }
                }{\end{list}}

\newcommand{\proof}[1]{{\noindent \bf Proof.~~}#1 \hspace{1cm} 
\medskip$\Box$\hfill}

\input{psfig}

\begin{document}
\thispagestyle{empty}
\maketitle
\SBIMSMark{1996/2}{March 1996}{}

\section{Introduction}
\label{teichintro}

In order to prove the renormalization conjecture for infinitely
renormalizable real polynomials of bounded combinatorics, Sullivan in
\cite{su3} introduced a space of analytic maps where the
renormalization operator is defined: the space of polynomial-like maps
of degree two and bounded combinatorics modulo holomorphic conjugacies
(see \cite{ms}). In this space it is possible to define a distance,
$\teichdist$ called the Teichm\"{u}ller distance. This distance
measures how far two polynomial-like maps are from being
holomorphically conjugate (see Definition~\ref{teichdistdef}).

It is obvious from the definition that the Teichm\"{u}ller distance is
a pseudo-distance. It is not obvious that this pseudo-distance is
actually a distance. To prove this is a distance it is necessary to
show that if two polynomial-like maps $f$ and $g$ are such that
$\teichdist (f,g)=0$ then they are holomorphically conjugate (this
can be viewed as a rigidity problem).
Sullivan showed in \cite{su3} that for real polynomials with connected
Julia set this is true. He makes use of external classes of
polynomial-like maps (as defined in \cite{dh}) to reduce the original
rigidity problem to a rigidity problem of expanding maps of the
circle, previously studied in \cite{ssu}. The last result concerning
expanding maps of the circle depends on the theory of Thermodynamical
formalism.

In this work we will show that the Teichm\"{u}ller distance for all elements
of a certain class of generalized polynomial-like maps (the class of {\em \ochs
generalized polynomial-like maps}, see Definition~\ref{pollikedef}) is
actually a distance, as in the case Sullivan studied.  This class contains
several important classes of generalized polynomial-like maps, namely: Yoccoz,
Lyubich, Sullivan and Fibonacci.

The initial motivation for this work was to carry on Sullivan's
renormalization theory for higher degree Fibonacci maps; the
renormalization scheme and $\apri$ bounds for these maps were given in
\cite{l} and \cite{lm}. The 
only place which needed adjustment after that was exactly the issue of
the Teichm\"{u}ller distance.

In our proof we can not use external arguments (like external
classes). Instead we use hyperbolic sets inside the Julia sets of our
maps. Those hyperbolic sets will allow us to use our main analytic
tool, namely Sullivan's rigidity Theorem for non-linear analytic
hyperbolic systems stated in Section~\ref{sulrigsection}.

Let us denote by $m$ the probability measure of maximal entropy for
the system $f:
\julf \rightarrow \julf$. In \cite{l0} Lyubich constructed
a maximal entropy measure $m$ for $f:
\julf \rightarrow \julf$ for any rational function $f$. Zdunik
classified in \cite{z} exactly when $\hausdim (m) = \hausdim ( \julf
)$. We show that the strict inequality holds if $f$ is \och, except
for Chebyshev polynomials.  This result is a particular case of
Zdunik's result if we consider $f$ as a polynomial.  It is however an
extension of Zdunik's result if $f$ is a generalized polynomial-like
map. The proof follows from the non-existence of invariant affine
structure proved in Section~\ref{non-affine}.

The structure of the paper is as follows: in
Section~\ref{teichstatement} we give the definitions necessary to
state the main Theorem. We also give the precise definition of the
class of polynomial-like maps that we will be working with. In
Section~\ref{thermodyn} we introduce notations and results concerning
Thermodynamical formalism which will be used later. As it was
mentioned before, in Section~\ref{sulrigsection} we give the statement
of Sullivan's rigidity Theorem for non-linear analytic hyperbolic systems. In
Section~\ref{hypset} we show that we can apply our main Theorem to
several classes of polynomial-like maps, namely: Sullivan, Yoccoz and
Lyubich polynomial-like maps and Fibonacci generalized polynomial-like
maps. We also in this section construct special hyperbolic sets inside
$\julf$. In Section~\ref{non-affine} we prove that the hypothesis of
 Sullivan's rigidity Theorem is satisfied
for the hyperbolic sets constructed in
Section~\ref{hypset}.  In Section~\ref{teichproof} we present the
proof of the Main Theorem. We finish the paper with
Section~\ref{non-affineconseq} where we derive some other consequences
from the result in Section~\ref{non-affine}.

\vspace{1cm}

\noindent {\bf Acknowledgment.} I would like to thank Misha Lyubich for 
many hours of extremely beneficial mathematical conversations and
suggestions. I am greatful to John Milnor for reading previous
versions of this work and making useful remarks. I also thank Jan Kiwi
and Marco Martens for various valuable mathematical conversations. I
am grateful to CNPq-Brazil and the Department of Mathematics of the
State University of New York at Stony Brook for financial support.
This work was part of the author's PhD thesis presented in June 1995 at 
the State University of New York at Stony Brook under the direction of
M.~Lyubich.

\section{Statement of the result}
\label{teichstatement}

\begin{DefinitioN}
\label{pollikedef}
Let $U$ and $U_{i}$ be open topological discs, $i=0, 1, ...,n$.
Suppose that $cl(U_{i}) \subset U$ and $U_{i} \bigcap U_{j} =
\emptyset$ if $i$ is different than $j$. A generalized polynomial-like
map is a map $\pollikef$ such that the restriction $f|U_{i}$ is a
branched covering of degree $d_{i}, d_{i} \geq 1$.
\end{DefinitioN}

We will not use the above Definition in full generality. From now on,
all generalized polynomial-like maps in this work will have just one
critical point. We will fix our notation as follows: $f|U_{0}$ is a
branched covering of degree $d$ onto $U$ (with zero being the only critical
point) and $f|U_{i}$ is an isomorphism onto $U$, if $i=1,...,n$.

The {\em filled in Julia set} of $f$, denoted by $\finjulf$, is
defined, as usually, as $\finjulf = \bigcap f^{-n} (\bigcup U_{i})$. The
{\em Julia set} of $f$, denoted by $\julf$, is defined as $\julf =
\partial (\finjulf )$. Douady and Hubbard introduced in \cite{dh} the
notion of a polynomial-like map. Their definition coincides with the
previous one when the domain of $f$ has just one component (the
critical one).  They also showed that a polynomial-like map of degree
$d$ is hybrid conjugate to a polynomial of the same degree (in some
neighborhoods of the Julia sets of the polynomial and the
polynomial-like maps). The above definition was given in \cite{l1}. It
was also showed that a generalized polynomial-like map is hybrid
conjugate to a polynomial (generally of higher degree but with only
one non-escaping critical point). We should keep in mind that
a polynomial map is a particular case of a generalized polynomial-like
map.

\begin{DefinitioN}
\label{classdef}
A generalized polynomial-like map $f$ is said to be \ochs if 
for any neighborhood of the critical point, the set of
points of $J(f)$ which avoid this neighborhood under the dynamics is
hyperbolic. We also ask $f$ to have its critical point in its Julia
set.
\end{DefinitioN}

There are several important examples of generalized polynomial-like
maps which are \och. Some example
are the following (see Lemma~\ref{hyperblemma}): Sullivan polynomials
(see \cite{ms} and
\cite{su3}), Yoccoz polynomials (see \cite{m}), Lyubich polynomials
(see \cite{l}) and their respective analog classes of polynomial-like
maps. Fibonacci generalized polynomial-like maps of even degree
(see \cite{lm}) are \ochs too.
All those classes just mentioned can
be put together inside one bigger class: the class of generalized
polynomial like maps which have $\apri$ bounds on infinitely many
generalized renormalization levels as described in
\cite{l} (which include the usual renormalization levels). 

Notice that there exist examples of polynomials which are not \och. 
That would be the case for $f$ having a 
neutral fixed point inside its Julia set (a Cremer polynomial, for
example).

\begin{DefinitioN}
\label{chebdef}
We say that a generalized polynomial-like map is Chebyshev if its
domain is connected (i.e., it is a polynomial-like map in the sense of
Douady and Hubbard) and the second iteration of its critical point is
a fixed point.
\end{DefinitioN}

\begin{DefinitioN}
\label{germdef}
We say that two generalized polynomial-like maps are in the same
conformal class if they are holomorphically conjugate in some
neighborhoods of their Julia sets. The conformal class of $f$ will be
denoted by $[f]$.
\end{DefinitioN}

\begin{DefinitioN}
\label{teichdistdef}
Let $[f]$ and $[g]$ be, as defined above, two conformal classes of
generalized polynomial-like maps. Let $h_{0}$ be a homeomorphism
conjugating $f$ and $g$ inside their respective Julia sets. Suppose
that there exist $U$ and $V$ neighborhoods of the Julia sets of
$f$ and $g$ and $h: U
\rightarrow V$ a conjugacy between $f$ and $g$. Assume that
$h$ is quasi-conformal with dilatation $K_{h}$ and that it is an
extension of $h_{0}$. Then we define the Teichm\"{u}ller distance
between $[f]$ and $[g]$ as $\teichdist ([f],[g]) = \inf_{h} \log
K_{h}$, where the infimum is taken over all conjugacies $h$
as described.
\end{DefinitioN}

Notice that $\teichdist ([f],[g]) \geq 0$ and $\teichdist([f],[g]) \leq
\teichdist([f],[t]) + \teichdist([t],[g])$, where $f$, $g$ and $t$ are
polynomial-like maps. In order to say that ``$\teichdist$''
is a distance we need to show that if $\teichdist([f],[g])=0$ then
$[f]=[g]$.
We prove the following Theorem:

\begin{thm}
\label{teichdistthm}
Let f and g be two generalized polynomial-like maps which are \och,
but not Chebyshev. Suppose that $\teichdist ([f],[g]) =0$. Then f and
g are conformally conjugate on a neighborhood of their Julia sets.
\end{thm}

We would like to point out that the above result is true even if $f$ 
is Chebyshev. That would not follow from our proof. It follows from 
the work of Sullivan \cite{su3}.

\section{Elements of thermodynamical formalism}
\label{thermodyn}

We refer the reader to \cite{b} for a detailed introduction to the
classical theory of thermodynamical formalism. See also \cite{pu} for
a more modern exposition of the subject. The goal of this Section is
to introduce notations and classical facts.

\begin{DefinitioN}
\label{hypsetdef}
Let $f$ be any conformal map. In what follows by a hyperbolic or
expanding set for $f$ we understand, as usual, a closed set $X$ such that
$f(X) \subset X$ and $|D(f^{n})(x)| \geq c \kappa^{n}$, for any $x$ in
$X$ and for $n \geq 0$, where $c>0$ and $\kappa>1$.
\end{DefinitioN}

\begin{DefinitioN}
\label{analrepdef}
We say that $f:X \rightarrow X$ is
topologically transitive if there exists a dense orbit in $X$.
\end{DefinitioN}

Suppose that the system $f:X
\rightarrow X$ is hyperbolic.
Then transitivity is equivalent to the following: for every non-empty
set $V \subset X$ open in $X$, there exists $n \geq 0$ such that
$\bigcup_{k \leq n} f^{k}(V)=X$. That is due to the existence of
Markov partition for $f: X \rightarrow X$.

Throughout this section the system $f:X \rightarrow X$ will be
conformal (i.e., $f$ will be defined and conformal in a neighborhood
of $X$) and hyperbolic.  If $\phi : X
\rightarrow {\rm R}$ is a H\"{o}lder continuous function, we say that
the probability measure $\mu_{\phi}$ is a {\em Gibbs measure}
associated to $\phi$ if:

\[ \sup_{\nu} \{ \entr_{\nu}(f) + \int_{X} \phi d \nu \} = 
\entr_{\mu_{\phi}}(f) + \int_{X} \phi d \mu_{\phi}  \]

\noindent where $\entr_{\nu}(f)$ is the entropy of $f$ with respect to 
the measure $\nu$ and the supremum is taken over all ergodic
probability measures $\nu$ of the system $f:X \rightarrow X$. In this
context we call $\phi$ a {\em potential function}. The {\em pressure
v} of the potential $\phi$ is denoted $ \preas ( \phi )$ and defined as
$\preas ( \phi ) = \sup_{\nu} \{ \entr _{\nu}(f)+ \int_{X} \phi d \nu
\}$, the supremum is taken over all ergodic probability measures.

The following Theorem assures us the existence of Gibbs measures.

\begin{theorem}[Ruelle-Sinai]
\label{bowen}
Given $f:X \rightarrow X$ hyperbolic and a H\"{o}lder continuous
potential $\phi : X \rightarrow {\rm R}$, there exists a unique Gibbs
measure $\mu_{\phi}$ associated to this potential.
\end{theorem}

One needs to know when two potentials generate the same Gibbs
measure. We have the following definition and Theorem to take care of that:

\begin{DefinitioN}
\label{cohomdef}
We say that two real valued functions $\phi , \psi : X \rightarrow {\rm
R}$ are cohomologous (with respect to the system $f:X \rightarrow X$)
if there exists a continuous function $s:X \rightarrow {\rm R}$ such that
$\phi (x) = \psi(x) + s(f(x))-s(x)$.
\end{DefinitioN}

\begin{theorem}[Livshitz]
\label{livthm}
Given $f:X \rightarrow X$ hyperbolic and two H\"{o}lder continuous
functions $\phi , \psi : X \rightarrow {\rm R}$, the following are
equivalent:

\begin{romanEnum}
\item $\mu_{\phi} = \mu_{\psi}$;
\item $\phi - \psi$ is cohomologous to a constant;
\item For any periodic point $x$ of $f:X \rightarrow X$ we have:

\[ \sum_{i=0}^{n-1} \phi (f^{i} (x) ) - \sum_{i=0}^{n-1} \psi (f^{i} (x) ) 
= n( \preas (\phi)-\preas (\psi)) \]
where $n$ is the period of $x$.
\end{romanEnum}
\end{theorem}

Of special interest is the one parameter family of potential functions
given by $\phi_{t}(x) =- t \log ( |Df (x)| )$. Notice that by the
definition of hyperbolic set, the functions $\phi_{t}$ are H\"{o}lder
continuous. One can study the pressure function $\preas (t) = \preas
( \phi _{t})$. Here are some properties of this function:

\begin{romanEnum}

\item $\preas(t)$ is a convex function;
\item $\preas(t)$ is a decreasing function;
\item $\preas (t)$ has only one zero exactly at $t= \hausdim (X)$;
\item $\preas (0) = \entr (f) =$ topological entropy of $f: X \rightarrow X$.
\end{romanEnum}

One can show that if $f: X \rightarrow X$ is hyperbolic, as we are
assuming, then the Hausdorff measure of $X$ is finite and positive.
That is because one can show that the Gibbs measure associated to the
potential given by $\phi_{\hausdim (X)} =- \hausdim(X) \cdot \log
(|Df|)$ is equivalent to the Hausdorff measure of $X$. Notice that
$\preas(\phi_{\hausdim (X)})=0$.  The Gibbs measure $\mu_{\phi_{0} -
\preas (\phi_{0})}$ associated to the potential $\phi_{0} -
\preas (\phi_{0}) \equiv -\preas ( \phi_{0}) = const$ is the measure of 
maximal entropy for the system $f: X \rightarrow X$. 
Instead of denoting this
measure by $\mu_{\phi_{0} - \preas (\phi_{0})}$ we will simply write
$\mu_{const}$.

Let us denote $m = \mu_{const}$ and $\nu = 
\mu_{\phi_{\hausdim
(X)}}$. The following is a consequence of the previous paragraph
and Theorem~\ref{livthm}.

\begin{corollary}
\label{confmaxentrmeascor}
Let $f:X \rightarrow X$ be hyperbolic.
The measures $m$ and $\nu$ are equal if and only if there exists
a number $\lambda$ such that for any periodic point $x$ of $f: X \rightarrow X$
we have $|Df^{n}(x)| = \lambda ^{n}$, where $n$ is the period of $x$.
\end{corollary}

\section{Sullivan's rigidity Theorem}
\label{sulrigsection}

We refer the reader to \cite{su2} and \cite{pu} for the proofs of the
results in this Section. In this section, the system $f:X \rightarrow X$
is assumed to be conformal (in a neighborhood of $X$) and hyperbolic.

\begin{DefinitioN}
\label{affinestrucdef}
An invariant affine structure for the system $f:X \rightarrow X$ is an atlas $
\{ (\sigma_{i}, U_{i}) \}_{i \in I}$ such that $\sigma_{i}: U_{i} \rightarrow
{\rm C}$ is a conformal injection for each $i$ where $X \subset
\bigcup_{i} U_{i}$ and all the maps $\sigma_{i} \sigma_{s}^{-1}$ and
$\sigma_{i} f \sigma_{s}^{-1}$ are affine (whenever they are defined).
\end{DefinitioN}

\begin{lemma}[Sullivan]
\label{sulcohomlemma}
Let $f:X \rightarrow X$ be a conformal transitive hyperbolic system. 
The potential
$\log (|Df|)$ is cohomologous to a locally constant function if and
only if $f:X \rightarrow X$ admits an invariant affine structure.
\end{lemma}

We call $f:X \rightarrow X$ a {\em non-linear system} if it does not
admit an invariant affine structure. Let $g:Y \rightarrow Y$ be
another system and let $h:X \rightarrow Y$ be a conjugacy between $f$
and $g$. Then we say that $h$ {\em preserves multipliers} if for every
$f$-periodic point of period $n$ we have $|Df^{n}(x)|
=|Dg^{n}(h(x))|$.

\begin{theorem}[Sullivan]
\label{sulrigthm}
Let $f:X \rightarrow X$ and $g:Y \rightarrow Y$ be two conformal non-linear
transitive hyperbolic systems. Suppose that $f$ and $g$ are conjugate by a
homeomorphism $h:X \rightarrow Y$ 
preserving multipliers. Then $h$ can be extended
to an analytic isomorphism from a neighborhood of $X$ onto a
neighborhood of $Y$.
\end{theorem}

\section{Hyperbolic sets inside the Julia set}
\label{hypset}

Let $\pollikef$ be any \ochs generalized polynomial-like map. Let $N$
be any neighborhood of the critical point.  We define:

\[ A_{N}= \{ z \in \julf : f^{j}(z) \notin N, \: \forall j \geq 0 \}\]
 
\noindent Notice that the set $A_{N}$ is forward $f$-invariant.
As $f$ is \och , we know that $f: A_{N} \rightarrow A_{N}$
is hyperbolic.

The next Lemma will show us that several important examples of
generalized polynomial-like maps are \och.

\begin{lemma}
\label{hyperblemma}
The map $f$ restricted to $A_{N}$ is hyperbolic if $f$ is either a
Yoccoz polynomial-like map or a Lyubich polynomial-like map or a
Sullivan polynomial-like map or a Fibonacci generalized
polynomial-like map of even degree.
\end{lemma}

\proof{
This Lemma is true because we can construct puzzle pieces for the set
$A_{N}$ if $f$ belongs to one of the classes mentioned in the
statement of this Lemma. 
We will describe how to do that. 

Let $f$ be a polynomial-like map (the generalized polynomial-like
map case is identical, as we will see later). We can find
neighborhoods $N_{n}$ of the critical point such that their boundaries
are made out of pieces of equipotentials and external rays landing at
appropriate pre-images of periodic points of $f$. Moreover the
diameter of $N_{n}$ tends to zero as $n$ grows (reference for this
fact: \cite{h} or \cite{m} if $f$ is Yoccoz,
\cite{l} if $f$ is Lyubich, \cite{hj} if $f$ is Sullivan, \cite{ls} if
$f$ is Fibonacci).  Let us fix $N_{n_{0}}$ such that $ N_{n_{0}}
\subset N$. 

By construction, $\partial N_{n_{0}} \bigcap J(f)$ is a finite set of
pre-images of periodic points. So there exists $l_{0}$ such that
$\bigcup_{i=0}^{l_{0}} f^{i}(
\partial N_{n_{0}} \bigcap J(f))$ is a forward invariant set under $f$. 
The same happens with the set of external rays landing at points in
$\partial N_{n_{0}} \bigcap J(f)$. So if $R$ is the set $\partial
N_{n_{0}} \bigcap J(f)$ together with the rays landing at $\partial
N_{n_{0}} \bigcap J(f)$, then there exists $l_{1}$ such that
$I = \bigcup_{i=0}^{l_{1}} f^{i}(R)$ is an invariant set under the
dynamics of $f$.

Each connected component of the complex plane minus the set $I$,
intersecting $A_{N}$ and bounded by a fixed equipotential of $\julf$
is by definition a puzzle piece of level zero for $A_{N}$. So we have
a Markov partition of our set $A_{N}$.  We define the puzzle-pieces of
level $k$ as being the connected components of the $k^{th}$ pre-image
of the puzzle pieces of level zero that intersect $A_{N}$. We will
denote by $Y_{n}(z)$ the puzzle piece of level $n$ containing $z$.

Thickening the puzzle pieces of level zero as described in \cite{m} we will
obtain open topological disks $V_{0}, ..., V_{n}$ covering $A_{N}$. We will
make use of the Poincar\'{e} metric on $V_{i}$, $0 \leq i \leq n$. We do this
thickening procedure in such a way to end up with exactly two branches of
$f^{-1}$ on each $V_{i}$, $0 \leq i \leq n$.  Each one of those inverse
branches is a holomorphic map carrying some $V_{i}$ isomorphically into a
proper subset of some $V_{j}$. For each puzzle piece of level zero, it follows
that the inverse branches of $f$ shrink the Poincar\'{e} distance by a factor
$\lambda < 1$. That is because each puzzle piece of level zero is compactly
contained in some thickened puzzle piece. Now, as $Y_{n}(z)$ is a connected
component of the pre-image of some puzzle piece of level zero under $f^{n}$,
it is easy to see that $\diam ( Y^{n}(z))$ tends to zero, as $n$ grows. So,
$\bigcap_{n \geq 0} Y^{n}(z) = \{ z \}$, for any $z$ in $A_{N}$.

Now one can show the hyperbolicity of $f:A_{N} \rightarrow A_{N}$.
Let $z$ be any element of $A_{N}$. We can take $Y_{n}(z)$ with
arbitrarily small diameter. Then the map $f^{n} : Y_{n}(z) \rightarrow
Y_{0}(f^{n}(z))$ is an isomorphism. If $V_{i}$ contains
$Y_{0}(f^{n}(z))$, then the appropriate inverse branch $f^{-n} :
Y_{0}(f^{n}(z)) \rightarrow Y_{n}(z)$ can be extended to $V_{i}$.  As
$Y_{0}(f^{n}(z))$ is compactly contained inside $V_{i}$, by Koebe
Theorem we conclude that the map $f^{n} : Y_{n}(z) \rightarrow
Y_{0}(f^{n}(z))$ has bounded distortion (not depending on $n$).  So we
have a map that maps isomorphically a set of arbitrarily small
diameter to a set of definite diameter with bounded distortion.  Those
observations together with the fact that $A_{N}$ is compact yield
hyperbolicity.

The same type of argument can be carried out if $f$ is a Fibonacci
generalized polynomial-like map (the only case with disconnected
domain that we are considering). If this is the case, then the domain
of $f$ has more than one component. In that case the puzzle pieces of
level zero are the connected components of the domain of the map $f$.
The puzzle pieces of higher levels are the pre-images of the puzzle
pieces of level zero. It follows from \cite{lm} (in the degree two
case) and
\cite{ls} (in the even degree greater then two case) that the puzzle
pieces shrink to points. The rest of the proof is identical to the
previous case.
}

We would like to point out that the above proof works for any
quadratic infinitely renormalizable generalized polynomial-like map
with $\apri$ bounds. The proof would be exactly the same. We would
use results from
\cite{j} for the construction of the small neighborhoods of the
critical points containing just pre-periodic points and external rays
on its boundaries (see also Theorem I in \cite{l4}).

Let $f$ be any \ochs generalized polynomial-like map. We will now construct a
sequence of sets that we will call $B_{n}$.  As the sets $A_{N}$, the sets
$B_{n}$ will also be $f$-invariant. The systems $f:B_{n} \rightarrow B_{n}$
will be hyperbolic and transitive. This is the main reason why we will need
this new family of sets.

Let us select one periodic point $p_{i}$ from every 
non post-critical periodic orbit of $f$ inside $\julf$.

\begin{lemma}
\label{bnconstrlemma}
Let $f$ be an \ochs generalized polynomial-like map. Then 
it is possible to construct a sequence of sets $B_{n} \subset J(f)$
such that:

\begin{romanEnum}

\item Each set $B_{n}$ is $f$-forward invariant, compact and hyperbolic;
\item  For any  $i=1,2,...,n$, 
the set of pre-images of $p_{i}$  belonging to $B_n$ is dense 
in $B_{n}$;
\item $f|B_{n} : B_{n} \rightarrow B_{n}$ is topologically transitive;
\item $\bigcup_{n} B_{n}$ is dense inside $J(f)$.
\end{romanEnum}

\end{lemma}

\proof{
Let us start the construction of the sets $B_{n}$. 
Let the period of
$p_{i}$ be $n_{i}$. We denote the orbit of $p_{i}$ by $\orb (p_{i})$.

We define the set $B_{1}$ simply as being $\orb (p_{1} )$. We will now
define the set $B_{2}$. Let $M_{i}$ be a small neighborhood of $p_{i}$
such that $f^{-n_{i}}(M_{i}) \subset M_{i}$ $i=1,2$. Here
$f^{-n_{i}}(M_{i})$ stands for the connected component of the
pre-image of $M_{i}$ under $f^{-n_{i}}$ containing $p_{i}$.  There
exists a pre-image $y_{1}$ of $p_{1}$ (suppose that $f^{s_{1}} (
y_{1})= p_{1}$) inside $M_{2}$ and a pre-image $y_{2}$ of $p_{2}$
(suppose that $f^{s_{2}} ( y_{2})= p_{2}$) inside $M_{1}$. The orbit
$y_{i}, f(y_{i}), .... f^{s_{i}}(y_{i})= p_{i}$ will be called a
bridge from $\orb (p_{i})$ to $\orb (p_{j})$, for $i \neq j$. There
exists a small neighborhood $\widetilde{ M_{i}} \subset M_{i}$
containing $p_{i}$ such that $y_{i} \in f^{-s_{i}} (\widetilde{
M_{i}}) \subset M_{j}$, $i \neq j$. Notice that we are not using
post-critical periodic points in our construction. That implies that
all the pre-images we are taking are at a positive distance from the
critical point.

In what follows $i \in \{1,2 \}$. Consider the pull back of the set
$M_{i}$ along the periodic orbit $p_{i}, f(p_{i}), ...,
f^{n_{i}}(p_{i})=p_{i}$: $M_{i}=M_{i}^{0}, M_{i}^{-1}, ...,
M_{i}^{-n_{i}+1}, M_{i}^{-n_{i}}$. Here $M_{i}^{-k}=f^{-k}(M_{i})$ 
for $k= 0, 1, ..., n_{i}$. 
Consider also the pull back of the set $\widetilde{M_{i}}$ along the
orbit $y_{i}, f(y_{i}), ..., f^{s_{i}}(y_{i})=p_{i}$: $\widetilde{M_{i}} =
\widetilde{M_{i}}^{0},
\widetilde{M_{i}}^{-1}, ...\widetilde{M_{i}}^{-s_{i}}$. Here 
$\widetilde{M_{i}}^{-k} = f^{-k}(\widetilde{M_{i}})$ for $k=0,1, ...,
s_{i}$. We have the following collections of inverse branches of $f$:
the first collection is $f^{-1} :M_{i}^{-l} \rightarrow M_{i}^{-l-1}$, 
for $l=0, 1,..., n_{i}-1$
and the second is $f^{-1}:
\widetilde{M_{i}}^{-l} \rightarrow
\widetilde{M_{i}}^{-l-1}$ for $l=0, 1, ..., s_{i}-1$ (remember that 
$i \in \{1,2 \} $). 

The union of the two collections of inverse branches of $f$ described
in the previous paragraph will be called our ``selection'' of branches
of $f^{-1}$ for $B_{2}$ (notice that we are specifying the domain and
image of each one of the branches of $f^{-1}$ in our ``selection'').
Consider now the set of all possible pre-images of $p_{i}$, $i=1,2$
under composition of branches of $f^{-1}$ in our ``selection'' of
branches.  We define the set $B_{2}$ as being the closure of the set
of all such pre-images.

We define $B_{n}$ in a similar fashion: instead of letting $i$ in
last paragraphs to be just in $\{1, 2 \}$, we let $i$ to be in $\{
1,2,...,n \}$. For each $p_{i}$, $M_{i}$ is as before a small
neighborhood around $p_{i}$, $i=1,2, ..., n$. There exist $y_{i,j}$
pre-image of $p_{i}$ (suppose that $f^{s_{i,j}} ( y_{i,j})= p_{i}$)
contained inside $M_{j}$ (those points define the bridges
between any two distinct orbits). There exists a small neighborhood
$\widetilde{ M_{i,j}}$ of $p_{i}$ contained in $M_{i}$ such that $y_{i,j}
\in f^{-s_{i,j}} (\widetilde{ M_{i,j}}) \subset M_{j}$, $i \neq j$.
As for $B_{2}$ now we can define the suitable ``selection'' of branches of
$f^{-1}$ for $B_{n}$ in a similar way. Consider now the set of
all possible pre-images of $p_{i}$, $i=1,2, ...,n$, under composition
of branches of $f^{-1}$ in our ``selection'' of branches. We define
the set $B_{n}$ as being the closure of the set of all pre-images just
described. Notice that we can carry on our construction such that we
have $B_{n-1} \subset B_{n}$. This finishes the construction of the sets 
$B_{n}$. Let us prove their properties.

The invariance and compactness are true by construction.
Hyperbolicity follows because we are excluding from our construction
post-critical periodic points. That implies that the distance from the
set $B_{n}$ to the critical point of $f$ is strictly positive
(depending on $n$). Hyperbolicity follows now, as $f$ is \och.

Let us show the second property. It is clear that the pre-images of the
set $\{p_{1}, p_{2}, ..., p_{n} \}$ under the system $f: B_{n}
\rightarrow B_{n}$ is dense inside $B_{n}$. So, in order to show that
the pre-images of some $p_{i}$ are dense inside $B_{n}$, we just need
to show that given any $1 \leq j \leq n$, there exist pre-images of
$p_{i}$ arbitrarily close to $p_{j}$. That is true because there
exists a pre-image $y_{i,j}$ of $p_{i}$ inside $M_{j}$, the
neighborhood of $p_{j}$ used in the construction of $B_{n}$. If we
take all pre-images of $y_{i,j}$ along the periodic orbit of of
$p_{j}$ we will find pre-images of $p_{i}$ arbitrarily close to
$p_{j}$ (remember that all periodic points are repelling).

Let us show $(iii)$. By $(ii)$, inside any open set $V \neq \emptyset$,
there exists a pre-image of $p_{i}$, for each $i=1,2,..., n$. Then,
for some $m_{i}$, $f^{m_{i}}(V)$ is a neighborhood of $p_{i}$, for
each $p_{i}$. Let $x$ be any point in $B_{n}$.  By the construction of
$B_{n}$, there exist a $j$ and a positive $k$ such that $f^{-k}(x) \in
M_{j}$.  Pulling $f^{-k}(x)$ back along the orbit of $p_{j}$
sufficiently many times we will find a pre-image of $x$ inside
$f^{m_{j}}(V)$. That implies that for some positive $s$, $x \in
f^{s}(V)$. So we conclude that $ B_{n} \subset
\bigcup_{k \geq 0} f^{k}( V)$. 

The last property is obvious because $\bigcup_{n} B_{n}$ contains all
the periodic points inside $\julf$ with the exception of at most
finitely many (in the case that the critical point is pre-periodic).
}

\section{Non-existence of an affine structure}
\label{non-affine}

In this Section we will show that if $f$ is \och, but not Chebyshev, then $f:B_{n}
\rightarrow B_{n}$ does not admit an invariant affine structure
 for $n > n_{0}$, for some $n_{0}$ depending on $f$.

\begin{lemma}
\label{eqaffstrclemma}
Suppose that $f:B_{n} \rightarrow B_{n}$ and $f:B_{n+1}
\rightarrow B_{n+1}$ admit invariant affine structures. Then the 
invariant affine structure
in $B_{n+1}$ extends the invariant affine structure in $B_{n}$.

\end{lemma}

\proof{ We will start by taking $n=1$. Let $\{ (\phi_{i} ,V_{i} ) \}$ 
be a finite atlas of an invariant affine structure for $ f:B_{1}
\rightarrow B_{1}$ and let $\{ (\sigma_{j}, U_{j}) \}$ be a finite atlas of
an invariant affine structure for $f:B_{2} \rightarrow B_{2}$. We will
show that the collection $\{ (\sigma_{j}, U_{j}) \}
\bigcup \{ (\phi_{i} , V_{i} ) \}$ is an atlas of an invariant affine
structure for $f:B_{2} \rightarrow B_{2}$. Notice that the invariant
affine structure for $f:B_{1} \rightarrow B_{1}$ is unique, given by
the linearization coordinates of $p_{1}$.

Let us suppose that $V_{i} \bigcap U_{j} \neq \emptyset$. We will
check that the change of coordinates $\sigma_{j} (\phi_{i})^{-1}$ is
affine. Suppose that the closure of
 $V_{i} \bigcap U_{j} \bigcap B_{2}$ is empty. Then if
we shrink $U_{j}$ (to $U_{j} \setminus V_{i}$) 
we can act just as if $V_{i} \bigcap U_{j} = \emptyset$.
So we can assume that the closure of $B_{2} \bigcap V_{i} \bigcap
U_{j}$ is not empty. Let $x$ be an element belonging to this
intersection. We can assume for simplicity that $V_{i}$ is a chart in
$B_{1}$ containing $p_{1}$ and $U_{i}$ is a chart in $B_{2}$
containing $p_{1}$ (remember that $p_{i}$ is the enumeration of
periodic points used to construct the sets $B_{n}$ and that $B_{1}
\subset B_{2}$). As the affine structure for periodic orbits is
unique, we conclude that $(\sigma_{i}, U_{i})$ and $(\phi_{i} ,V_{i}
)$ are the same (up to an affine map) in a neighborhood of $p_{1}$.

We can pull $x$ back by $f^{n_{1}}$ along the (periodic) orbit of
$p_{1}$ until we find $y \in B_{2}$, a pre-image of $x$ under some
iterate of $f^{n_{1}}$. Notice that $y$ is in fact an element of $B_{2}$.
That is because the inverse branches of $f$ following the orbit of
$p_{1}$ are in the ``selection'' of inverse branches used to construct
$B_{2}$.  Because $(\phi_{i} ,V_{i} )$ is a linearization coordinate
around the periodic point $p_{1}$, we conclude that $\phi_{i}
f^{ln_{1}} (\phi_{i})^{-1}$ is affine from a neighborhood of
$\phi_{i}(y)$ to a neighborhood of $\phi_{i}(x)$. On the other hand,
as $y \in B_{2}
\bigcap U_{i}$ and $x=f^{ln_{1}}(y) \in B_{2} \bigcap U_{j}$, we
conclude that $\sigma_{j} f^{ln_{1}} (\sigma_{i})^{-1}$ is affine from
a neighborhood of $\sigma_{i}(y)$ to a neighborhood of
$\sigma_{j}(x)$. Keeping in mind that $\sigma_{i}$ is equal to
$\phi_{i}$ (up to an affine map), we get that the change of
coordinate $\sigma_{j} (\phi_{i})^{-1}$ is affine (see
Figure~\ref{t3}). From that follows trivially that any composition of
the form $\sigma_{j} f (\phi_{i})^{-1}$ and $\phi_{i} f
(\sigma_{j})^{-1}$ is affine, whenever they are defined. So we proved
the Lemma in the case $n=1$.

\begin{figure}[htb]
\centerline{\psfig{figure=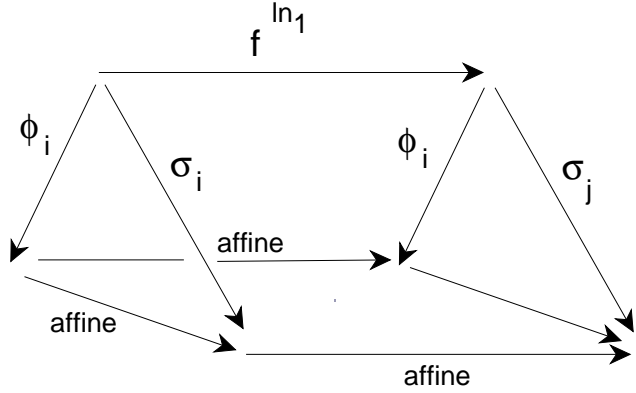,height=6cm}}
\caption{Commutative diagram of charts}
\label{t3}
\end{figure}

Now suppose that we have some invariant affine structure $\{ (\phi_{i}
,V_{i} ) \}$ in $B_{n}$ and $\{ (\sigma_{j}, U_{j}) \}$ in $B_{n+1}$.
We want to show that $\{ (\phi_{i} ,V_{i} ) \} \bigcup \{ (\sigma_{j},
U_{j}) \}$ is an invariant affine structure in $B_{n+1}$. Suppose that
$U_{j}$ intersects a chart of one of the periodic points $p_{1},
p_{2}, ..., p_{n}$ in $B_{n}$. Then the change of coordinates from
$U_{j}$ to one of those charts is affine (same as the proof for
$n=1$). Now let $U_{j}$ and $V_{i}$ be two arbitrary charts with
non-empty intersection. We can assume that there exists $x$ in the
closure of $B_{n+1} \bigcap U_{j} \bigcap V_{i}$ (otherwise we can
shrink $U_{j}$ to $U_{j} \setminus V_{i}$).  Let $V_{1}$ be the chart
around $p_{1}$ in the affine structure for $B_{n}$ and $U_{1}$ be the
chart around $p_{1}$ in the affine structure for $B_{n+1}$. Then
$V_{1} \bigcap U_{1}$ is a neighborhood of $p_{1}$. Inside $B_{n}$ we
can pull back $V_{i}$ until we find a pre-image (of $V_{i}$ with
respect to some iterate of $f: B_{n} \rightarrow B_{n}$) strictly
inside $U_{1}
\bigcap V_{1}$ (this is possible by property $(iii)$
 in Lemma~\ref{bnconstrlemma} and hyperbolicity).
So $f^{-l} (V_{i})
\subset V_{1} \bigcap U_{1} $. Then it is clear that $\phi_{i} f^{l}
(\phi_{1})^{-1}$ is affine in $f^{-l} (V_{i})$. On the other hand,
$\sigma_{j} f^{l} (\sigma_{1})^{-1}$ is affine in a subset of
$f^{-l} (V_{i})$ containing the pre-image $y$ of $x$ via $f^{-l}$
(remember that $x$ is the element in $B_{n+1}$ contained in $U_{j}
\bigcap V_{i}$). As $\sigma_{1}$ and $\phi_{1}$ are equal up to an
affine transformation in a neighborhood of $y$ (because both are
linearizing coordinates around a periodic point), we conclude that the
change of coordinates $\phi_{i} (\sigma_{j})^{-1}$ is affine (just
imagine an appropriate diagram similar to the one in Figure~\ref{t3}).
It is trivial to check that the affine structure defined by $\{
(\phi_{i} ,V_{i} ) \} \bigcup \{ (\sigma_{j}, U_{j}) \}$ is invariant
under $f$.  
}

We would like to point out that with exactly the same demonstration as
above we show that if there exists an invariant affine structure for
the system $f:B_{n} \rightarrow B_{n}$, then it is unique.

\begin{lemma}
\label{nonexist}
If $f$ is \ochs and is not Chebyshev, then there
is a positive number $n_{0}$ such that $f:B_{n} \rightarrow B_{n}$
does not admit an invariant affine structure if $n> n_{0}$ ($n_{0}$
depends on $f$ ).
\end{lemma}

\proof{
Suppose that $f:B_{n} \rightarrow B_{n}$ admits an affine structure,
for infinitely many $n$. Then all those structures coincide when
defined in common subsets by Lemma~\ref{eqaffstrclemma}. This implies that
 we can
define the set $X =
\bigcup_{n} B_{n}$ and an invariant affine 
structure for $f:X \rightarrow X$ (notice that $X$ is $f$-invariant
and dense inside $\julf$). Let us
denote the elements of the atlas defining such affine structure over
$X$ by $( \sigma_{i},U_{i})$.

There exists $n$ such that some element of $f^{-n} (0)$, say $y_{0}$,
belongs to some $U_{\beta}$ (here we need to have  our map $f$ not 
conjugate to $z^{d}$).  There exists $m$ such that some element of
$f^{-m}(f^{2}(0))$, say $y_{1}$, which is not a pre-image of the
critical point $0$ belongs to $U_{\alpha}$, for some $\alpha$ (notice
that this is not true if $f$ is a Chebyshev generalized 
polynomial-like map). We
can take $U_{\alpha}$ and $U_{\beta}$ small enough such that $f^{m}:
U_{\alpha}
\rightarrow f^{m} ( U_{\alpha}) =
\ U_{\alpha}' $ and $f^{n}: U_{\beta} \rightarrow f^{n} ( U_{\beta}) =  
U_{\beta} '$ are isomorphisms (see Figure~\ref{t1}). We can also
assume that $f^{2} (U_{\beta}') = U_{\alpha}'$. We can find $x \in X
\bigcap U_{\beta}'$ 
because $X$ is dense inside $\julf$. Notice that we need to have the 
critical point inside $J(f)$ in order to be able to construct 
$U_{\beta}'$ with small 
diameter intersecting $X$. Then $f^{2}(x) \in U_{\alpha}'$.
We can take charts from the atlas on $X$, say $(\sigma_{\gamma},
U_{\gamma}), U_{\gamma} \subset U_{\beta}'$ and $(\sigma_{\nu},
U_{\nu}), U_{\nu} \subset U_{\alpha'}$ containing $x$ and $f^{2} (x)$
respectively. Let $\sigma_{\beta}' = \sigma_{\beta} f^{-n}$ and
$\sigma_{\alpha}' = \sigma_{\alpha} f^{-m}$, where the inverse
branches $f^{-n}$ and $f^{-m}$ are defined according to our previous
discussion. Notice that $\sigma_{\beta}'$ and $\sigma_{\alpha}'$ are
isomorphisms onto their respective images. Let $A = \sigma_{\nu} f^{2}
\sigma_{\gamma}^{-1}$. The map $A$ is affine (because $A$ is the map $f^{2}$
viewed from the atlas over $X$).

\begin{figure}[htb]
\centerline{\psfig{figure=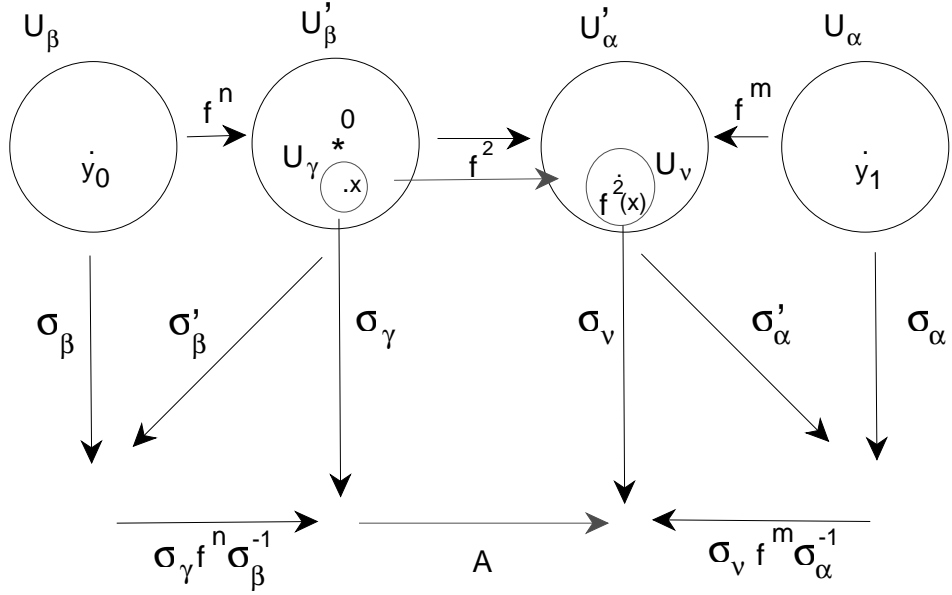,height=8.5cm}}
\caption{Commutative diagram}
\label{t1}
\end{figure}

Notice that

\[   \sigma_{\alpha}' f^{2} (\sigma_{\beta}')^{-1}  =  
(\sigma_{\nu} f^{m}
\sigma_{\alpha}^{-1})^{-1} A (\sigma_{\gamma} f^{n} \sigma_{\beta}^{-1}) \] 
when we restrict both sides of the
equation to the set $\sigma_{\beta}' ( U_{\gamma} )$.

The left-hand side of the above formula is a restriction of a degree
two branched covering, namely $\sigma_{\alpha} f^{-m} f^{2} f^{n} (
\sigma _{\beta})^{-1}$. If we restrict the right-hand side of our equation to
$\sigma_{\beta}' (U_{\gamma})$ we get an affine map. Contradiction! 
The Lemma is proved.  
}

\section{Proof of the Theorem}
\label{teichproof}

We will present in this Section the proof of Theorem~\ref{teichdistthm}. Let
$\pollikef$ and $\pollikeg$ be two \ochs generalized polynomial-like maps, but
not Chebyshev. Let us suppose that $\teichdist (f,g) = 0$. This implies that there
exists a homeomorphism $h: \julf \rightarrow \julg$ conjugating $f$ and $g$
which is extended by quasi-conformal maps of arbitrarily small
distortion. This implies that $h$ preserves multipliers.

We define the hyperbolic sets $X_{n} = B_{n} \subset \julf$ (as
introduced in Section~\ref{hypset}) and $Y_{n}=h(X_{n}) \subset
\julg$.

The systems $f:X_{n} \rightarrow X_{n}$ and $g:Y_{n} \rightarrow
Y_{n}$ do not admit invariant affine structures if $n$ is big (see
Lemma~\ref{nonexist}). In other words, $f:X_{n} \rightarrow X_{n}$
and $g:Y_{n} \rightarrow Y_{n}$ are non-linear systems, for $n$ big.
So by Theorem~\ref{sulrigthm} we
know that there exist open neighborhoods $O_{n}$ of $X_{n}$ and
$O_{n}'$ of $Y_{n}$ and holomorphic isomorphisms $H_{n}: O_{n}
\rightarrow O_{n}'$ extending $h_{n}$. We can assume that $O_{n}
\subset O_{n+1}$ and $O_{n}' \subset O_{n+1}'$. Notice that by
analytic continuation we have that $H_{n}=H_{n+1}$ on $O_{n}$.

We define two open sets, $O= \bigcup_{n} O_{n}$ and $O'=
\bigcup_{n} O_{n}'$. We can define $H: O \rightarrow O'$ by the
following: for any $z \in O$ there exists some $n$ such that $z \in
O_{n}$. Then we define $H(z)=H_{n}(z)$. The map $H$ is well defined.
The map $H$ is holomorphic because locally it coincides with $H_{n}$,
for some $n$.  It is also injective. The map $H$ conjugates $f|X$ and
$g|Y$, where we define $X= \bigcup X_{n}$ and $Y= \bigcup Y_{n}$. The
sets $X$ and $Y$ are dense subsets of $\julf$ and $\julg$,
respectively. So the conjugacy $H$ is defined in a open neighborhood
of a dense subset of $\julf$. Our goal is to extend $H$ to a
neighborhood of the whole Julia set.

Suppose that $z$ is a point in $\julf$ not belonging to $O$. If $z$ is
not the critical value, then there exists $n$ and an element $z_{-n}$
of $f^{-n}(z)$, such that $z_{-n} \in O$, and the iteration $f^{n}$
restricted to a small ball around $z_{-n}$ is injective. Consider the
holomorphic map defined in a small neighborhood $W$ of $z$ by $\phi =
g^{n} H f^{-n}$, where by $f^{-n}$ we understand the branch of
$f^{-n}$ that takes $z$ to $z_{-n}$. If $W$ is sufficiently small,
then $\phi$ is an isomorphism. It is clear that $\phi$ and $h$
coincide where both are defined. By analytic continuation, that means
that $\phi$ also coincides with $H$ where both are defined. In this
way we managed to extend $H$ to an open neighborhood of $\julf
\setminus \{ f(0) \}$. We will keep calling this extension $H$. If $z$
is the critical value, then instead of looking for pre-images of $z$
in order to repeat the previous reasoning, just look for the first
image of $z$.  Remember that now the second iterate of the critical
point belongs to the domain of $H$. We can define an isomorphism in a
small neighborhood of the critical value given by $\phi =g^{-1} H
f$. The same argument as before goes through to show that we have
extended $H$ to an open neighborhood of $\julf$.  This proves the
Theorem.

\section{Other consequences of the non-existence of affine structure}
\label{non-affineconseq}

According to Lemma~\ref{sulcohomlemma}, the non-existence of affine
structure for the system $f:A_{N} \rightarrow A_{N}$ is equivalent to
$\log(|Df)|$ not cohomologous to a locally constant function in
$A_{N}$. In particular, the non-existence of an affine structure
implies that $\log(|Df|)$ is not cohomologous to a constant function
inside $A_{N}$. This last observation together with
Theorem~\ref{livthm} implies the following:

\begin{corollary}
\label{multiplierlemma}
If $f$ is \och, then there is no $\lambda$ such
that for any $n$ and any $f$-periodic point $p$, $|Df^{n}(p)|= \lambda
^{n}$.
\end{corollary}

\proof{ By our previous comments, we conclude that if $\diam (N)$ 
is small, then there is no $\lambda$ such that $|Df^{n}(p)|= \lambda ^{n}$ 
for any $n$ and any
$f$-periodic point $p$ inside $A_{N}$. That
implies the Corollary.  }

If $\mu$ is a Borel probability measure in $\julf$, then we define the
Hausdorff dimension of $\mu$ as $\hausdim (\mu ) = \inf \hausdim (Y)$
where the infimum is taken over all sets $Y \subset \julf$ with $\mu
(Y) =1$.

Remember that the measure $m = \mu_{const}$ is the measure of maximal
entropy for the hyperbolic system $f:X \rightarrow X$. Zdunik proved 
in \cite{z} that for rational maps $\hausdim (m) = \hausdim ( \julf
)$ if and only if $f$ is $ z \mapsto z^{d}$ or a Chebyshev polynomial.  
The following is a
particular case of Zdunik's result if we consider $f$ as a
polynomial. It is however an extension of Zdunik's result if $f$ is a
generalized polynomial-like map:

\begin{corollary}
If $f$ is \ochs and $m$ is the measure of maximal
entropy for $f$, then $\hausdim (m) < \hausdim (J(f))$.
\end{corollary}

\proof{ It was shown in \cite{puz} (see Theorem 6) 
 that it is enough to check that $\log(|Df|)$ is not cohomologous to a
constant in $\julf$. By that we mean the following: there is no real
function $h$ which is equal $m$-a.e. to a continuous function in a
small neighborhood of any point in $\julf$ without the post-critical
set and $\log(|Df|) = c + h(f(x)) -h(x)$.

Suppose that $\log(|Df|)$ is cohomologous to a constant, in the sense
defined in the previous paragraph. Remember that the sets $B_{n}$ are
at a positive distance from the closure of the critical orbit. So we
would conclude that $\log(|Df|)$ is cohomologous to a constant (in the
sense of Definition~\ref{cohomdef}). Lemma~\ref{sulcohomlemma} and
Lemma~\ref{nonexist} imply that this is impossible.  }

\bibliography{todo}

\begin{thebibliography}{dMvS93}

\bibitem[Bow75]{b}
R.~Bowen.
\newblock {\em Equilibrium states and the ergodic theory of {A}nosov
  diffeomorphisms}, volume 470 of {\em Lecture Notes in Mathematics}.
\newblock Springer-Verlag, 1975.

\bibitem[DH85]{dh}
A.~Douady and J.~Hubbard.
\newblock On the dynamics of polynomial-like maps.
\newblock {\em Ann. Sc. \'{E}c. Norm. Sup}, 18, 287-343, 1985.

\bibitem[dMvS93]{ms}
W.~de~Melo and S.~van Strein.
\newblock {\em One dimensional dynamics}.
\newblock Springer-Verlag, 1993.

\bibitem[HJ]{hj}
J.~Hu and Y.~Jiang.
\newblock The {J}ulia set of the {F}eigenbaum quadratic polynomial is locally
  connected.
\newblock Preprint, 1993.

\bibitem[Hub]{h}
J.~Hubbard.
\newblock Local connectivity of {J}ulia sets and bifurcation loci: thee
  theorems of {J}.-{C}. {Y}occoz.
\newblock In {\em Topologycal methods in modern mathematics, A symposium in
  honor of {J}ohn {M}ilnor}. Publish or Perish, 467-511.

\bibitem[Jia]{j}
Y.~Jiang.
\newblock Infinitely renormalizable quadratic {J}ulia sets.
\newblock Preprint, 1993.

\bibitem[LM93]{lm}
M.~Lyubich and J.~Milnor.
\newblock The unimodal {F}ibonacci map.
\newblock {\em J. of the A. M. S.}, 6, 425-457, 1993.

\bibitem[LvS95]{ls}
G.~Levin and S.~van Strein.
\newblock Local connectivity of {J}ulia set of real polynomials.
\newblock {\em IMS@Stony Brook preprint series}, (1995/5), 1995.

\bibitem[Lyu83]{l0}
M.~Lyubich.
\newblock Entropy properties of rational endomorphisms of the {R}iemann sphere.
\newblock {\em Erg Th. and Dyn. Syst.}, 3, 351-385, 1983.

\bibitem[Lyu91]{l1}
M.~Lyubich.
\newblock On the {L}ebesgue measure of the {J}ulia set of a quadratic
  polynomial.
\newblock {\em IMS@Stony Brook preprint series}, (1991/10), 1991.

\bibitem[Lyu93]{l}
M.~Lyubich.
\newblock Geometry of quadratic polynomials: moduli, rigidity and local
  connectivity.
\newblock {\em IMS@Stony Brook preprint series}, (1993/9), 1993.

\bibitem[Lyu95]{l4}
M.~Lyubich.
\newblock Dynamics of quadratic polynomials {II}: Rigidity.
\newblock {\em IMS@Stony Brook preprint series}, (1995/14), 1995.

\bibitem[Mil91]{m}
J.~Milnor.
\newblock Local connectivity of {J}ulia sets: expository lectures.
\newblock {\em IMS@Stony Brook preprint series}, (1991/10), 1991.

\bibitem[PU]{pu}
F.~Przyticki and M.~Urba\'{n}ski.
\newblock {\em Monograph in preparation}.

\bibitem[PUZ89]{puz}
F.~Przyticki, M.~Urba\'{n}sk, and A.~Zdunik.
\newblock Harmonic, {G}ibbs and {H}ausdorff measures on repellers for
  holomorphic maps {I}.
\newblock {\em Ann. of Math.}, 130, 1-40, 1989.

\bibitem[SS85]{ssu}
M.~Shub and D.~Sullivan.
\newblock Expanding endomorphism of the circle revisited.
\newblock {\em Erg. Th. and Dyn. Syst.}, 5, 285-289, 1985.

\bibitem[Sul86]{su2}
D.~Sullivan.
\newblock {\em Quasiconformal homeomorphisms in dynamics, topology and
  geometry}.
\newblock Proc I. C. M., Berkeley. 1986.

\bibitem[Sul92]{su3}
D.~Sullivan.
\newblock Bounds, quadratic differentials and renormalization conjecture.
\newblock In {\em A. M. S. centenial publication 2: Mathematics into the
  Twenty-first century}. A. M. S., 1992.

\bibitem[Zdu90]{z}
A.~Zdunik.
\newblock Parabolic orbifolds and the {H}ausdorff dimension of maximal entropy
  measures for rational maps.
\newblock {\em Inv. Math.}, 996, 27-649, 1990.

\end{thebibliography}
\bibliographystyle{alpha}

\end{document}